\documentclass[11pt]{article}

\usepackage{a4wide}
\usepackage{amssymb}
\usepackage{amsfonts}
\usepackage{amsmath}
\input xy
\xyoption{arrow}
\xyoption{matrix}

\date{}

\newtheorem{proposition}{Proposition}[section]
\newtheorem{theorem}[proposition]{Theorem}
\newtheorem{lemma}[proposition]{Lemma}

\newtheorem{corollary}[proposition]{Corollary}

\def\GK{{\rm  GK}\,}
\def\Kdim{{\rm K.dim }\,}

\def\Hom{{\rm Hom}}
\def\der{\partial }

\def\nFM0{{\nu }_{F,M_0}}
\def\nFN0{{\nu }_{F,N_0}}
\def\nGN0{{\nu }_{G,N_0}}

\def\N0{ {\bf N}_0 }

\def\t{\otimes}

\def\ra{\rightarrow}

\def\lra{\leftrightarrow}
\def\Xpm{X^{\pm }}

\def\s{\sigma}
\def\Z{\mathbb{Z}}

\def\l1{{\lambda}_1}

\def\a{\alpha}
\def\a0{ {\alpha }_0}
\def\a1{ {\alpha }_1}

\def\l{\lambda}

\def\nFGM0{{\nu }_{F,G,M_0}}

\def\nFN0{{\nu}_{F,N_0}}

\def\sm{{\sigma}^m}

\def\sm1{{\sigma}^{-1}}

\def\smtp1{{\sigma}^{-t+1}}

\def\S1{S^{-1}}

\def\Xpm1{X^{\pm 1}_1}

\def\sPM1{{\sigma }^{\pm 1}}
\def\sMP1{{\sigma }^{\mp 1 }}

\def\d{\delta}

\def\di{{\rm d.ind}}

\def\L{\Lambda}

\def\OO{{\cal O}}

\def\CD{{\cal D}}

\def\Ytm1{Y^{t-1}}
\def\Yim1{Y^{i-1}}

\def\k{{\bf k}}

\def\i{{\bf i}}

\def\Der{{\rm Der }}
\def\ad{{\rm ad }}

\def\ker{ {\rm ker } }

\def\D{ \Delta }

\def\SL2Z{ {\rm SL}_2({\bf Z}) }

\def\Gp1{ G^{1 , 1 } }
\def\P11{ P^{-1 , 1 } }
\def\Pp1{ P^{1 , 1 } }

\def\nCLsr{{}^\nu\kern-2pt {\cal L}^{\sigma , \rho  }}
\def\nP{{}^\nu \kern-2pt P}
\def\nL{{}^\nu\kern-2pt L}
\def\nLL{{}^\nu\kern-2pt \Lambda}
\def\nPsr{{}^\nu\kern-2pt P^{\sigma , \rho  }}
\def\nLsr{{}^\nu\kern-2pt L^{\sigma , \rho  }}
\def\nuCL{{}^\nu\kern-2pt  {\cal L}}
\def\nCLsr{{}^\nu\kern-2pt {\cal L}^{\sigma , \rho  }}
\def\nCL1m{{}^\nu\kern-2pt {\cal L}^{-1 , 1  }}
\def\x1nu{x^\frac{1}{\nu}}
\def\xm1nu{x^{-\frac{1}{\nu}}}

\def\trdeg{{\rm tr.deg}}

\def\ra{\rightarrow }

\def\nAM0{{\nu }_{{\cal A},M_0}}
\def\nAN0{{\nu }_{{\cal A},N_0}}

\def\Kdim{ {\rm Kdim } }

\def\End{ {\rm End }}
\def\Der{ {\rm Der }}

\def\det{ {\rm det }}
\def\ad{ {\rm ad }}

\def\bp{\overline{p}}

\def\ga{\mathfrak{a}}

\def\bJ{\overline{J}}

\def\derij{\partial_{{\bf i}, {\bf j}}}

\def\j{{\bf j}}

\def\Spec{{\rm Spec}}

\def\Hom{{\rm Hom}}
\def\II{{\bf I}}
\def\JJ{{\bf J}}

\def\HS{{\rm HS}}
\def\di!{\frac{\der^i}{i!}}
\def\dj!{\frac{\der^j}{j!}}
\def\dik!{\frac{\der^k_i}{k!}}

\def\HS{{\rm HS}}
\def\di!{\frac{\der^i}{i!}}
\def\dik!{\frac{\der^k_i}{k!}}

\def\N{\mathbb{N}}

\def\0{\overline{0}}
\def\1{\overline{1}}

\def\Ln1{\L_{n,\overline{1}}}

\def\a1{a_{\overline{1}}}

\def\S{\Sigma}

\def\vn1{\overrightarrow{n-1}}

\def\id{{\rm id}}

\begin{document}

\author{V. V. \  Bavula 
}

\title{Generators and defining relations for
 ring of differential operators on  smooth affine algebraic
variety in prime characteristic}

\maketitle
\begin{abstract}
For the ring of differential operators $\CD (\OO (X))$ on a smooth
affine algebraic variety $X$ over a perfect field of
characteristic $p>0$, a set of algebra generators and a  set of
defining relations are found explicitly.
 A  finite set of generators and a finite set of
 defining relations are given explicitly for the module
$\Der_K(\OO (X))$ of derivations on the algebra $\OO (X)$ of
regular functions on the variety $X$. For an arbitrary irreducible
affine algebraic variety $X$,  it is proved that each term $\CD
(\OO (X))_i$ of the  order filtration $\CD (\OO (X))=\cup_{i\geq
0}\CD (\OO (X))_i$ is a finitely generated left $\OO (X)$-module.
The same results are true for  ring  of differential operators on
 regular algebra of essentially finite type.

\end{abstract}


\section{Introduction}

In prime characteristic, differential operators and their modules
are a more difficult  and less developed area of Mathematics then
in characteristic zero. The main difficulty is that in prime
characteristic algebras of differential operators are not finitely
generated, not left or right Noetherian, and contain a lot of
nilpotent elements. As a result methods of affine Algebraic
Geometry are not applicable (at least in the way they are in
characteristic zero), and this is a principal unavoidable problem.

Key ingredients of the theory of (algebraic) $\CD $-modules in
 characteristic zero are the {Gelfand-Kirillov dimension,
 multiplicity, Hilbert polynomial, the inequality of Bernstein, and
 holonomic modules}. In prime characteristic, straightforward
 generalizations of these either do not exist or give `wrong'
 answers (as in the case of the Gelfand-Kirillov dimension:  $\GK
 (\CD (P_n))=n$ in prime characteristic rather than $2n$ as it
 `should' be and it is in characteristic zero where $P_n$ is a polynomial algebra in $n$ variables).

 In 70's and 80's, for rings of differential operators in prime characteristic
  natural questions were posed (see, for example, questions 1-4 in
 \cite{Smith85LNM}) [some of them are still open] that can be
 summarized as {\em to find generalizations of the mentioned
 concepts and results (that results in `good theory'} expectation
 of which was/is high, see, the remark of ${\rm Bj\ddot{o}rk}$ in
 \cite{Smith85LNM}). One of the questions in the paper of P. Smith  \cite{Smith85LNM} is
 {\em to give a definition of holonomic module in prime
 characteristic}. In characteristic zero, holonomic modules have
 remarkable {\em homological} properties based on which
   Mebkhout and Narvaez-Macarro \cite{Meb-Nar-MacLNM90} gave a
 definition of holonomic module. Another approach (based on the
 {\em Cartier Lemma}) was taken by Bogvad \cite{BogvadJA95} who defined, so-called,
 {\em filtration holonomic} modules. This one is more close to the
 original idea of holonomicity in characteristic zero.
 Note that the two mentioned concepts of holonomicity in prime
 characteristic appeared {\em before} analogues of the
 Gelfand-Kirillov dimension and the inequality of Bernstein have
 been recently  found, \cite{holmodp}. In \cite{holmodp},  analogues of the Gelfand-Kirillov dimension,
multiplicity, the inequality of Bernstein, and holonomic modules
are found in prime characteristic (using very different methods
and ideas from characteristic zero) and classical properties of
holonomic modules were proved.  In prime characteristic,
$\CD$-modules were studied by Haastert \cite{Haastert87},
Alvarez-Montaner,  Blickle and  Lyubeznik \cite{Alv-Blickle-Lyub},
Bezrukavnikov,  Mirkovic, and  Rumynin
\cite{Bezrukavnikov-Mirkovic-Rumynin-AnnMAth}; and they were used
in study of (local) cohomology by Huneke and Sharp
\cite{Hun-Sharp93}, Kashiwara and Lauritzen
\cite{Kashiwara-Lauritzen2002}, Andersen and
 Kaneda \cite{Andersen-Kaneda2000},   and Lyubeznik
\cite{LyubCrelle97}. $\CD$-modules  were applied to the theory of
tight closure by K. Smith \cite{K.Smith95} and to the ring of
invariants by K. Smith and van den Bergh \cite{KSm-vdB}.

$\noindent $

In this paper, module means a {\em left} module, $\N := \{ 0,1, 2,
\ldots \}$ is the set of natural numbers. The following notation
will remain fixed throughout the paper (if it is not stated
otherwise):
\begin{itemize}
\item $K$ is a {\bf perfect} field of characteristic $p>0$ (not
necessarily algebraically closed);  \item  $P_n:=K[x_1, \ldots ,
x_n]$ is a polynomial algebra over $K$; \item
$\der_1:=\frac{\der}{\der x_1}, \ldots , \der_n:=\frac{\der}{\der
x_n}\in \Der_K(P_n)$;  \item $I:=\sum_{i=1}^m P_nf_i$ is a {\bf
prime} but {\bf not} a maximal ideal of the polynomial algebra
$P_n$ with a set of generators $f_1, \ldots , f_m$;  \item  an
algebra $A:=S^{-1} (P_n/I)$ (the localization of the algebra
$P_n/I$ at a multiplicatively closed subset $S$ of $P_n/I$) which
is a domain with the field of fractions $Q:={\rm Frac}(A)$, i.e.
$A$ is an arbitrary algebra of essentially finite type which is a
domain; \item $\Der_K(A)$ is the left $A$-module of
$K$-derivations of the algebra $A$; \item $\HS_K(A)$ is the set of
higher derivations (Hasse-Schmidt derivations) of the algebra $A$;
 \item  $\CD
(A)$ is the ring of $K$-linear differential operators on the
algebra $A$. The action of a differential operator $\d \in \CD
(A)$ on an element $a\in A$ is denoted either by $\d (a)$ or $\d
*a$;
 \item the homomorphism $\pi :P_n\ra A$, $p\mapsto \bp $, to
make notation simpler we sometime write $x_i$ for $\overline{x}_i$
(if it does not lead to confusion);  \item  the {\bf Jacobi}
$m\times n$ matrices
 $J=(\frac{\der f_i}{\der x_j})\in M_{m,n}(P_n)$ and
 $\bJ =(\overline{\frac{\der f_i}{\der
x_j}})\in M_{m,n}(A)\subseteq M_{m,n}(Q)$;   $r:={\rm rk}_Q(\bJ )$
is the  rank of the Jacobi matrix $\bJ$ over the field $Q$; \item
$\ga_r$ is the {\bf Jacobian ideal} of the algebra $A$ which is
(by definition) generated by all the $r\times r$ minors of the
Jacobi matrix $\bJ$ ($A$ is {\em regular} iff $\ga_r=A$, it is the
{\bf Jacobian criterion of regularity}, \cite{Eisenbook}, 16.20);
\end{itemize}

For $\i =(i_1, \ldots , i_r)$ such that $1\leq i_1<\cdots <i_r\leq
m$ and $\j =(j_1, \ldots , j_r)$ such that $1\leq j_1<\cdots
<j_r\leq n$, $\D
 (\i , \j )$ denotes the corresponding minor of the Jacobi matrix
$\bJ =(\bJ_{ij})$, that is $\det (\bJ_{i_\nu, j_\mu})$, $\nu , \mu
=1, \ldots, r$;  and the element $\i$ (resp. $\j $) is called {\bf
non-singular} if $\D (\i , \j')\neq 0$ (resp. $\D (\i', \j )\neq
0$) for some $\j'$ (resp. $\i'$). We denote by $\II_r$ (resp.
$\JJ_r$) the set of all the non-singular $r$-tuples $\i$ (resp.
$\j $).

Since $r$ is the rank of the Jacobi matrix $\bJ$, it is easy to
show that $\D (\i , \j )\neq 0$ iff $\i\in \II_r$ and $\j\in
\JJ_r$ (Lemma \ref{IIr1}). Denote by $\JJ_{r+1}$ the set of all
$(r+1)$-tuples $\j =(j_1, \ldots , j_{r+1})$ such that $1\leq
j_1<\cdots <j_{r+1}\leq n$ and when deleting  some element, say
$j_\nu$, we have a  non-singular $r$-tuple $(j_1, \ldots
,\widehat{j_\nu},\ldots , j_{r+1})\in \JJ_r$ where hat over a
symbol here and everywhere means that the symbol is omitted from
the list. The set $\JJ_{r+1}$ is called the {\bf critical set} and
any element of it is called a {\bf critical singular}
$(r+1)$-tuple.

Let us describe the main results of the paper.  The next theorem
gives a set of generators and a
 set of defining relations for the left $A$-module
$\Der_K(A)$ when $A$ is a  regular algebra.
\begin{theorem}\label{9bFeb05}
Let the algebra $A$ be a regular algebra. Then the left $A$-module
$\Der_K(A)$ is generated by derivations $\der_{\i , \j }$, $\i \in
\II_r$, $\j \in \JJ_{r+1}$, where
\begin{eqnarray*}
 \derij  =  \der_{i_1, \ldots , i_r; j_1, \ldots , j_{r+1}}:= \det
 \begin{pmatrix}
  \overline{\frac{\der f_{i_1}}{\der x_{j_1}}} & \cdots &  \overline{\frac{\der f_{i_1}}{\der
  x_{j_{r+1}}}}\\
  \vdots & \vdots & \vdots \\
  \overline{\frac{\der f_{i_r}}{\der x_{j_1}}} & \cdots &  \overline{\frac{\der f_{i_r}}{\der
  x_{j_{r+1}}}}\\
  \der_{j_1}& \cdots & \der_{j_{r+1}}\\
\end{pmatrix}
\end{eqnarray*}

that satisfy the following defining relations (as a left
$A$-module): 
\begin{equation}\label{Derel}
\D (\i , \j )\der_{\i', \j'}=\sum_{l=1}^s(-1)^{r+1+\nu_l}\D (\i';
j_1', \ldots , \widehat{j_{\nu_l}'}, \ldots , j_{r+1}')\der_{\i;\j
, j_{\nu_l}'}
\end{equation}
for all $\i, \i'\in \II_r$, $\j=(j_1, \ldots , j_r)\in \JJ_r$, and
$\j'=(j_1', \ldots , j_{r+1}')\in \JJ_{r+1}$ where $\{ j_{\nu_1}',
\ldots , j_{\nu_s}'\}=\{ j_1', \ldots , j_{r+1}'\}\backslash \{
j_1, \ldots , j_r\}$.
\end{theorem}

For elements $\i = (i_1, \ldots , i_r)\in \II_r$ and $\j = (j_1,
\ldots , j_r)\in \JJ_r$, Theorem \ref{pCDAd=gd} gives commuting
iterative higher derivations $\{ \d_{\i ; \j ,
j_\nu}^{[k]}\}_{k\geq 0}$, $\nu = r+1, \ldots , n$ of the
localization $A_\D$ of the algebra $A$ at the powers of the
element $\D := \D (\i , \j )$ where $ \{ j_{r+1}, \ldots , j_n\}=
\{ 1, \ldots , n \}\backslash \{ j_1, \ldots , j_r\}$.  These
iterative higher derivations can be found explicitly, i.e. their
actions on the generators $x_i$ can be found explicitly (Theorem
\ref{rpCDAd=gd} and Theorem \ref{12Jul8}). Note that each higher
derivation is uniquely determined by its action on algebra
generators since the corresponding automorphism does.

Choose a function $\N \ra \N$, $ k\mapsto n(k)$ such that
$n(0)=0$,  
\begin{equation}\label{nk1}
n(k+l) \geq n(k) + n(l) \;\; {\rm for \; all}\;\; k,l\in \N ;
\end{equation}
\begin{equation}\label{dk1}
d_{\i ; \j , j_\nu}^{[k]}:= \D^{n(k)} \d_{\i ; \j ,
j_\nu}^{[k]}\in \CD (A), \;\; k\geq 0,
\end{equation}
for all $\i\in \II_r$, $\j \in \JJ_r$ and $\nu = r+1, \ldots , n$.
Note that $d_{\i ; \j , j_\nu}^{[0]}=\id_A$, the identity map on
$A$. The function $n(k)$ can be found explicitly (see Section
\ref{Reldefops}). Moreover, any fast growing function
 satisfies the conditions (\ref{nk1}) and  (\ref{dk1}).

The next result gives a set of generators and a set of {\em
defining} relations for the $K$-algebra $\CD (A)$ of differential
operators on $A$ (it is well known and is not difficult to show
that the algebra $\CD (A)$ is not finitely generated and does not
satisfy finitely many defining relations).

\begin{theorem}\label{21Jul8}
Let the algebra $A$ be a regular algebra. Then the ring of
differential operators $\CD (A)$ on $A$ is a simple algebra
generated over $K$ by the algebra $A$ and the elements $$d_{\i ;
\j , j_\nu}^{[k]}:=\D^{n(k)} \d_{\i ; \j , j_\nu}^{[k]},\;\; k\geq
1,\;\;  \i\in \II_r,\;\;  \j \in \JJ_r,\;\;  \nu = r+1, \ldots ,
n$$
  where
 $ \{ j_{r+1}, \ldots , j_n\}=
\{ 1, \ldots , n \}\backslash \{ j_1, \ldots , j_r\}$ and  $\D :=
\D (\i , \j )$. These elements satisfy the defining relations
(R1)--(R5) over $A$: for all $\nu, \mu = r+1, \ldots , n$ and
natural numbers $k,l\geq 1$,

(R1) $$\D^{n(l)}d_{\i ; \j , j_\nu}^{[k]}\sum_{s=0}^l a_{\i ; \j ,
j_\mu}^{[l]}(0,s) d_{\i ; \j , j_\mu}^{[s]}\, \D^{n(k)+k}=
\D^{n(k)}d_{\i ; \j , j_\mu}^{[l]}\sum_{t=0}^k a_{\i ; \j ,
j_\nu}^{[k]}(0,t) d_{\i ; \j , j_\nu}^{[t]}\, \D^{n(l)+l}$$ for
some elements $a_{\i ; \j , j_\mu}^{[l]}(0,s), a_{\i ; \j ,
j_\nu}^{[k]}(0,t)\in A$;

(R2) $$ \D^{n(k+l)-n(k)}d_{\i ; \j , j_\nu}^{[k]}\sum_{s=0}^l
a_{\i ; \j , j_\nu}^{[l]}(0,s) d_{\i ; \j , j_\nu}^{[s]} =
{k+l\choose l} d_{\i ; \j , j_\nu}^{[k+l]}\D^{n(l)+l}; $$

(R3) $$[d_{\i ; \j , j_\nu}^{[k]}, x_{j_\mu }]= \d_{j_\nu , j_\mu
} \D^{n(k) - n(k-1)}d_{\i ; \j , j_\nu}^{[k-1]} $$ where
$\d_{j_\nu , j_\mu } $ is the Kronecker delta;

(R4) $$[d_{\i ; \j , j_\nu}^{[k]}, x_{j_s
}]=\sum_{t=1}^k\D^{n(k)-n(t)-n(k-t)}d_{\i ; \j ,
j_\nu}^{[t]}(x_{j_s })d_{\i ; \j , j_\nu}^{[k-t]}, \;\; s=1,
\ldots , r;$$

(R5) $$ \D^{m(l)}d_{\i' ; \j' ,
j'_\s}^{[l]}\D^{\S_{n-r}}=\D'^{n(l)}\sum_{|\k |\leq l} c_{l, \k }
\prod_{\rho = 1}^{n-r}\sum_{t_{r+\rho } =0}^{k_{r+\rho }} a_{\i ;
\j , j_\nu }^{[k_{r+\rho }]}(\S_{\rho -1}, t_{r+\rho }) d_{\i ; \j
, j_\nu}^{ [t_{r+\rho }] }$$ for all $\i'\in \II_r$, $\j'\in
\JJ_r$, $j'_\s \in \{ 1, \ldots , n\} \backslash \{ j_1', \ldots ,
j_r'\}$ and some  $c_{l, \k }= c_{l, \k } (\i'; \j',
 j'_\s ;  \i , \j )\in A$,  and $m(l)\in \N$ where $\D':= \D (\i',
\j')$, $\S_0:=0$, and
$$ \S_\rho := \sum_{\nu =1}^\rho (n(k_{r+\nu} )+k_{r+\nu}), \;\;\;  \rho
\geq 1, \;\;\;  \k := (k_{r+1}, \ldots , k_n) \in \N^{n-r}, |\k
|:= k_{r+1}+\cdots +k_n. $$
\end{theorem}

{\it Remark}. All the elements of the algebra $A$ in the defining
relations (R1)--(R5) are found explicitly, see (\ref{R4Ds2}),
(\ref{clk1}) and (\ref{blk2}).


$\noindent $

 When  $S=\{ 1\}$, the algebra $A= P_n/I$ is the algebra of regular functions on the irreducible
 affine algebraic variety $X=\Spec (A)$, therefore we have the
 explicit algebra generators for the ring of differential
 operators $\CD (X)= \CD (A)$ on an arbitrary  smooth irreducible affine
 algebraic variety $X$. Any regular affine algebra $A'$ is a
 finite direct product of regular affine domains,
 $A'=\prod_{i=1}^s A_i$.  {\em  Since $\CD (A')\simeq \prod_{i=1}^s
 \CD (A_i)$, Theorem \ref{21Jul8} gives  algebra generators and defining relations for the ring of
 differential operators on arbitrary smooth affine algebraic
 variety. Since $\Der_K(A')\simeq \bigoplus_{i=1}^s\Der_K(A_i)$,
 Theorem \ref{9bFeb05} gives generators and defining relations for the left
 $A'$-module of derivations} $\Der_K(A')$ of the algebra $A$.

$\noindent $

In characteristic zero, analogues of Theorems \ref{9bFeb05} and
\ref{21Jul8} were proved in \cite{gendif}, Theorems 1.1 and 1.2
respectively. Theorem 1.1, \cite{gendif} coincides with Theorem
\ref{9bFeb05}, but Theorem 1.2, \cite{gendif} is much more simpler
than Theorem \ref{21Jul8} (the main difference is that the
defining relations in Theorem 1.2, \cite{gendif} are of {\em first
order} in derivations).

\begin{theorem}\label{9Feb05}
({\rm Theorem
1.2,} \cite{gendif}) Let the algebra $A$ be a regular algebra over
a field $K$ of characteristic zero. Then the ring of differential
operators $\CD (A)$ on $A$ is generated over $K$ by the algebra
$A$ and the derivations $\der_{\i , \j }$, $\i \in \II_r$, $\j \in
\JJ_{r+1}$ that satisfy the defining relations
(\ref{Derel}) and 
\begin{equation}\label{1Derel}
\der_{\i , \j }\overline{x}_k=\overline{x}_k\der_{\i , \j
}+\der_{\i , \j }(\overline{x}_k), \; \i \in \II_r, \; \j \in
\JJ_{r+1}, \; k=1, \ldots , n.
\end{equation}

\end{theorem}

{\it Definition}. The {\bf higher derivation algebra} $\D (A)$ is
a subalgebra of the ring of differential operators $\CD (A)$
generated by the algebra $A$ and the higher  derivations
$\HS_K(A)$ of the algebra $A$.

\begin{theorem}\label{23Jul8}
(Criterion of regularity of the algebra $A$  via $\D (A)$) The
following statements are equivalent.
\begin{enumerate}
\item $A$ is a regular algebra. \item $\D (A)$ is a simple
algebra.\item $A$ is a simple $\D (A)$-module.
\end{enumerate}
\end{theorem}

In characteristic zero, the same criterion was proved in
\cite{MR}, 15.3.8 where $\D (A)$ is {\bf the derivation algebra},
it is a subalgebra  of the ring of differential operators $\CD
(A)$ generated by the algebra $A$ and the set  $\Der_K(A)$ of all
the $K$-derivations of the algebra $A$.

If the field $K$ has characteristic  zero and  the algebra $A$ is
{\em regular} then the ring of differential operators $\CD (A)$ is
a finitely generated Noetherian algebra. If $A$ is {\em not}
regular then, in general, the algebra $\CD (A)$ need not be  a
finitely generated algebra nor a left or right Noetherian algebra,
\cite{BGGDiffcone72}; the algebra $\CD (A)$ can be finitely
generated and right Noetherian yet not left Noetherian,
\cite{SmStafDifopcurve} (so, in characteristic zero for a
non-regular algebra $A$ the ring $\CD (A)$ behaves similarly as
the ring $\CD (A)$ for a regular algebra $A$ in prime
characteristic). Though, a kind of finiteness still holds for a
{\em singular} algebra $A$ in characteristic zero.

\begin{theorem}\label{CDAifg}
{\rm (Theorem 1.5,
\cite{gendif})}  Let $K$ be a field of  characteristic  zero and
$\CD (A)=\cup_{i\geq 0}\CD (A)_i$ be the order filtration of $\CD
(A)$. Then, for each $i\geq 0$, $\CD (A)_i$ is a finitely
generated left $A$-module.
\end{theorem}

In Section \ref{opsing}, we prove the same result in prime
characteristic.

\begin{theorem}\label{pCDAifg}
 Let $K$ be a perfect  field of  characteristic
$p>0$   and $\CD (A)=\cup_{i\geq 0}\CD (A)_i$ be the order
filtration of $\CD (A)$. Then, for each $i\geq 0$, $\CD (A)_i$ is
a finitely generated left $A$-module.
\end{theorem}


\section{Generators and defining relations for  ring of differential operators
on regular algebra of essentially finite
type}\label{Reldefops}

In this section, Theorems \ref{9bFeb05}, \ref{21Jul8} and
\ref{23Jul8}  are proved.

$\noindent $

Let $B$ be a commutative $K$-algebra. The ring of ($K$-linear)
{\bf differential operators} $\CD (B)$ on $B$ is defined as a
union of $B$-modules  $\CD (B)=\cup_{i=0}^\infty \,\CD (B)_i$
where $\CD (B)_0=\End_R(B)\simeq B$, ($(x\mapsto bx)\lra b$),
$$ \CD (B)_i=\{ u\in \End_K(B):\, [r,u]:=ru-ur\in \CD (B)_{i-1}\; {\rm for \; each \; }\; r\in B\}.$$
 The set of $B$-modules $\{ \CD (B)_i\}_{i\geq 0}$ is the {\em order filtration} for
the algebra $\CD (B)$:
$$\CD(B)_0\subseteq   \CD (B)_1\subseteq \cdots \subseteq
\CD (B)_i\subseteq \cdots\;\; {\rm and}\;\; \CD (B)_i\CD
(B)_j\subseteq \CD (B)_{i+j}, \;\; i,j\geq 0.$$

\begin{lemma}\label{IIr1}
$\i \in \II_r$ and $\j\in \JJ_r$ $\Leftrightarrow$ $\D (\i , \j
)\neq 0$.
\end{lemma}

{\it Proof}. Repeat the proof of the same statement in
characteristic zero, Lemma 2.1, \cite{gendif}. $\Box$

$\noindent $

Theorem \ref{p25July04} describes the derivations of an {\em
arbitrary} (not necessarily regular) domain of essentially finite
type.

\begin{theorem}\label{p25July04}
 Let   $\i
=(i_1, \ldots , i_r)\in \II_r$,  $\j =(j_1, \ldots , j_r)\in
\JJ_r$, and $\{ 1,\ldots , n \}\backslash \{ j_1, \ldots , j_r\}=$
$\{ j_{r+1}, \ldots , j_n\} $. Then   $\Der_K(A)=\{ \D (\i , \j
)^{-1} \sum_{k=r+1}^n a_{j_k}\der_{i_1, \ldots , i_r; j_1, \ldots
, j_r,j_k}\, | $ where the elements $ a_{j_{r+1}}, \ldots ,
a_{j_n}\in A$ satisfy the following system of inclusions:
$$ \sum_{k=r+1}^n\D (\i ; j_1, \ldots , j_{\nu -1},
j_k, j_{\nu +1}, \ldots , j_r)a_{j_k}\in A\D (\i , \j ), \;\; \nu
=1, \ldots , r\}. $$
\end{theorem}

{\it Proof}. Repeat the proof of the same statement in
characteristic zero, Theorem 2.12.(1), \cite{gendif}. $\Box$

$\noindent $

 Let us recall
basic facts about higher derivations. For more detail the reader
is referred to \cite{Ma}, Sec. 27.

A sequence $\d =(1:={\rm id}_A, \d_1, \d_2, \ldots )$ of
$K$-linear maps from a $K$-algebra $A$ to itself is called a {\em
higher derivation} (or a {\em Hasse-Schmidt derivation}) over $K$
from $A$ to $A$ if, for each $k\geq 0$,  
\begin{equation}\label{dkxyhi}
\d_k(xy)=\sum_{i+j=k}\d_i(x)\d_j(y)\; \;\; {\rm for\; all}\;\;
x,y\in A.
\end{equation}
 These conditions are equivalent to
saying that the map
$$e:A\ra A[[t]], \;\; x\mapsto \sum_{i\geq
0}\d_i(x)t^i,$$
 is a $K$-algebra homomorphism where $A[[t]]$ is a
ring of power series with coefficients from $A$. Let $\HS_K(A)$ be
the set of higher derivations on $A$. In general, a higher
derivation $\d = (\d_i)$ is not determined by the derivation
$\d_1$.

Let  $\d = (\d_i)\in \HS_K(A)$. By (\ref{dkxyhi}), $\d_1\in
\Der_K(A)$ and 
\begin{equation}\label{diCDAi}
\d_i\in \CD (A)_i, \;\; i\geq 0,
\end{equation}
since $\d_ix-x\d_i= \sum_{j=0}^{i-1}\d_{i-j}(x)\d_j$ for all $x\in
A$ and the result follows by induction on $i$.

A higher derivation $\d =(\d_i)\in \HS_K(A)$ is called {\em
iterative} if $\d_i\d_j={i+j\choose i}\d_{i+j}$ for all $i,j\geq
0$. Then a  direct computation shows that 
\begin{equation}\label{dip=0}
\d_{i}^p=0\;\;\; {\rm  for \; all\;} \;\;  i\geq 1,
\end{equation}
$\d_{i}^p=\d_{i}\cdots \d_{i}= {2i\choose i}{3i\choose 2i}\cdots
{pi\choose (p-1)i}\d_{pi}=0\d_{pi}=0$. For $i=1$, we have
$\d_1^p=0$.

{\bf The higher derivations $(1, \frac{\der_i}{1!},
\frac{\der^2_i}{2!}, \ldots )\in \HS_K(P_n)$, $i=1,\ldots , n$}.
The $K$-algebra homomorphism $P_n\ra P_n[[t]]$, $f(x_1, \ldots ,
x_n)\mapsto f(x_1, \ldots , x_{i-1}, x_i+t, x_{i+1}, \ldots ,
x_n)=\sum_{i\geq 0}\dik! (f)t^k$, gives the higher derivation
 $(1, \frac{\der_i}{1!},
\frac{\der^2_i}{2!}, \ldots )\in \HS_K(P_n)$. If ${\rm char}
(K)=0$ then $\dik!$ means $(k!)^{-1} \der^k_i$, but if ${\rm char}
(K)=p>0$ then 
\begin{equation}\label{nbinijchp}
\dik! (x_j^l)=\d_{ij}{l\choose k}x_j^{l-k}
\end{equation}
for all $k\geq 1$, $l\geq 0$, and $1\leq i,j\leq n$, where
$\d_{ij}$ is the Kronecker delta.

The action of the higher derivations $\der_i^{[k]} :
=\frac{\der_i^k}{k!} $ on the polynomial algebra $P_n=K\t_\Z \Z [
x_1, \ldots , x_n]$ should be understood as the action of the
element $1\t_\Z\frac{\der_i^k}{k!}$. The higher derivations $\{
\der_i^{[k]} \}_{k\geq 0}$ are iterative and they commute,
$\der_i^{[k]}\der_j^{[l]}= \der_j^{[l]}\der_i^{[k]}$. For each
element $\alpha = (\alpha_i) \in \mathbb{N}^n$, let
$\der^{[\alpha]}:= \prod_{i=1}^n \der_i^{[\alpha_i]}$.

\begin{theorem}\label{CDAffl}
Let the algebra $A$ be a regular algebra. Then $\CD (A)\ra
\prod_{\i\in \II_r , \j\in \JJ_r}\CD (A)_{\D (\i , \j )}$ is a
left and right faithfully flat extension of algebras where $\CD
(A)_{\D (\i , \j )}$ is the localization of the algebra $\CD (A)$
at the powers of the element $\D (\i , \j )$.
\end{theorem}

{\it Proof}. The algebra $A$ is regular, so $A=\ga_r=(\D (\i , \j
))=(\D (\i , \j ))_{\i\in \II_r, \j\in \JJ_r}$, hence the ideal of
$A$ generated by any power of the elements $\{ \D_{\i , \j}\, | \,
\i\in \II_r, \j\in \JJ_r\} $ is also equal to $A$. The extension
is a flat monomorphism. Suppose that the extension is not, say
left faithful, then there exists a proper left ideal, say $L$, of
$\CD(A)$ such that $\prod_{\i\in \II_r , \j\in \JJ_r}\CD (A)_{\D
(\i , \j )}\t_{\CD (A)}(\CD (A)/L)=0$, equivalently, there exists
a sufficiently large natural number $k$ such that $\D (\i , \j
)^k\in L$ for all $\i \in \II_r, \j\in \JJ_r $. Since $A=(\D (\i ,
\j )^k)_{\i \in \II_r, \j\in \JJ_r}\subseteq L$, we must have
$L=\CD (A)$, a contradiction. $\Box $

Let $R$ be a (not necessarily commutative)  algebra over a field
$K$, and let $\d $ be a
 $K$-derivation of the algebra $R$. For any elements $a,b\in R$
 and a natural number $n$, an easy induction argument gives the
 Leibniz formula
 $$ \d^n(ab)=\sum_{i=0}^n\, {n\choose i}\d^i(a)\d^{n-i}(b).$$
 It follows that the kernel $C(\d , R):=\ker \, \d $ of $\d $ is a
 subalgebra (of constants for $\d $) of $R$ (since $\d (ab)=\d
 (a)b+a\d (b)=0$ for $a,b\in C(\d , R)$), and the union of the
 vector spaces $N(\d ,R)=\cup_{i\geq 0}\, N(\d , i,R)$ is a positively
 filtered algebra (so-called, the {\em nil-algebra} of $\d $) where $N(\d
 , i,R):=\{ a\in R\, | \, \d^{i+1}(a)=0\}$,  that is
 $$N(\d , i,R)N(\d , j,R)\subseteq N(\d , i+j,R), \;\; {\rm for \;\;
 all}\;\; i,j\geq 0.$$
 Clearly, $N(\d , 0,R)= C(\d , R)$ and $N(\d , R):=\{ a\in R \, | \ \d^n (a)=0$
  for some natural $n\}$.

A $K$-derivation $\d $ of
 the algebra $R$ is a {\em locally nilpotent } derivation if for
 each element $a\in R$ there exists a natural number $n$ such
 that $\d^n(a)=0$. A $K$-derivation $\d $ is locally nilpotent iff
 $R=N(\d , R)$. A derivation of $R$ of the type $\ad (r) :x\mapsto
 [r,x]:=rx-xr$ is called an {\em inner} derivation of $R$ where
 $r\in R$.

\begin{lemma}\label{pdx=1}
(Lemma 2.1, \cite{JPAA-08invcharp})  Let $R$ be an algebra over an
arbitrary field $K$, $\d $ be a $K$-derivation of $R$ such that
$\d (x_i)=x_{i-1}$, $i\geq 0$ for some elements $x_i\in R$ such
that $x_{-1}=0$ and $x_1=1$. Then $N(\d , R)=\bigoplus_{i\geq 0}
Cx_i=\bigoplus_{i\geq 0} x_iC$ where  $C:=\ker \, \d $, and $N(\d
, i, R)=\bigoplus_{j=0}^i\, Cx_j=\bigoplus_{j=0}^i\, x_jC$ for all
$i\geq 0$.
\end{lemma}

%

As the first application of Lemma \ref{pdx=1}, we find generators
and defining relations for the algebra of differential operators
with polynomial coefficients. The results are known but the proof
is new.
 We have included the proof since it is short and similar patterns
 will appear later in the proofs of similar results in the general
 situation (Theorems \ref{pCDAd=gd} and \ref{21Jul8}).

\begin{corollary}\label{c1pdx=1}

\begin{enumerate}
\item The algebra $\CD (P_n)$ of differential operators with
polynomial coefficients $P_n$ is generated (as an abstract
$K$-algebra) by the elements $x_i$, $\der_i^{[k]}$, $i=1, \ldots ,
n$ and $k\geq 1$, that satisfy the following defining relations:
for all $i,j=1, \ldots , n$ and $k,l\geq 0$,
$$ [x_i, x_j]=[\der_i^{[k]},\der_j^{[l]}]=0,\;\;\;
\der_i^{[k]}\der_i^{[l]}={k+l\choose k}\der_i^{[k+l]} \;\;\;
[\der_i^{[k]}, x_j]=\d_{ij}\der_i^{[k-1]}, $$ where $\d_{ij}$ is
the Kronecker delta. \item $ \CD (P_n)=\bigoplus_{\alpha, \beta
\in \mathbb{N}^n}P_n x^\alpha \der^{[\beta]} =\bigoplus_{\alpha,
\beta \in \mathbb{N}^n}P_n \der^{[\beta]} x^\alpha $ and $ \CD
(P_n)_i=\bigoplus_{|\alpha |+| \beta |\leq i}P_n x^\alpha
 \der^{[\beta]} =\bigoplus_{|\alpha |+|
\beta |\leq i}P_n  \der^{[\beta]} x^\alpha $ for $i\geq 0$. \item
The algebra $\CD (P_n)$ is a central  simple algebra generated by
the polynomial algebra $P_n$ and its higher derivations
$HS_K(P_n)$. 
\item The map $\CD (P_n)\ra \CD
(P_n)^o$, $x_i\mapsto x_i$, $\der_i^{[k]}\mapsto (-1)^k
\der_i^{[k]}$, is a $K$ algebra isomorphism where $\CD (P_n)^o$ is
the opposite algebra. So, the algebra $\CD (P_n)$ is self-dual.
\end{enumerate}
\end{corollary}

{\it Proof}. It is obvious that the elements $x_i$, $\der_i^{[k]}$
where  $i=1, \ldots , n$ and $k\geq 1$ satisfy the given
relations. Then applying Lemma \ref{pdx=1} several times to the
algebra $E:={\rm End}_K(P_n)$ and the set of commuting inner
derivations $\ad \, x_1, \ldots , \ad \, x_n$, we obtain the
algebra
$$ N:=N(\ad \, x_1, \ldots , \ad \, x_n; E):=\cap_{i=1}^nN(\ad \,
x_i, E)=\bigoplus_{\alpha \in
\mathbb{N}^n}C\der^{[\alpha]}=\bigoplus_{\alpha \in
\mathbb{N}^n}\der^{[\alpha]}C$$ where $\alpha =(\alpha_1, \ldots ,
\alpha_n)$, $\der^{[\alpha]}:= \prod_{i=1}^n \der_i^{[\alpha_i]}$,
and $C:=\cap_{i=1}^n\ker (\ad \, x_i)={\rm End}_{P_n}(P_n)\simeq
P_n$. Therefore, $$ N=\bigoplus_{\alpha, \beta \in
\mathbb{N}^n}P_n x^\alpha \der^{[\beta]} =\bigoplus_{\alpha, \beta
\in \mathbb{N}^n}P_n \der^{[\beta]} x^\alpha .$$ It follows that
$N\subseteq \CD (P_n)$, the inverse inclusion follows at once from
the definitions of the algebra $N$ and $\CD (P_n)$. Therefore,
$N=\CD (P_n)$. It follows that $$ \CD (P_n)_i=\bigoplus_{|\alpha
|+| \beta |\leq i}P_n x^\alpha \der^{[\beta]} =\bigoplus_{|\alpha
|+| \beta |\leq i}P_n \der^{[\beta]} x^\alpha $$ for $i\geq 0$.
This proves statement 2 and statement 3 apart from simplicity of
the algebra $\CD (P_n)$.

To prove simplicity of the algebra $\CD (P_n)$, let $a=\sum
a_{\alpha \beta} x^\alpha \der^{[\beta]} $ be a nonzero element of
the algebra $\CD (P_n)$ where $a_{\alpha \beta} \in K$. We have to
show that the ideal $(a)$ generated by the element $a$ is equal to
the algebra $\CD (P_n)$. To prove this we use induction on the
degree $d:=\deg (a)=\max \{ |\alpha | +|\beta |\, | \, a_{\alpha
\beta} \neq 0\}$ of the element $a$. The case $d=0$ is obvious.
Suppose that the result is true for all nonzero elements of degree
$<d$. If there exists a coefficient $a_{\alpha \beta}\neq 0$ for
some $\beta \neq 0$, i.e. $\beta_i\neq 0$ for some $i$, then
applying the inner derivation $\ad \, x_i$ to the element $a$ we
have a nonzero element $[x_i, a]$ of degree $<d$, then induction
gives the result.

Now, we are in a situation where $a_{\alpha \beta} =0$ for all
$\beta \neq 0$, that is $a\in P_n$ is a polynomial of degree
$d>0$. Then there exists a variable, say $x_i$, such that
$\deg_{x_i}(a)=m>0$ (the degree in $x_i$). Then applying the inner
derivation $\ad (\der_i^{[m]})$ to the element $a$ we have a
nonzero element of degree $<d$, and induction finishes the proof
of simplicity of the algebra $\CD (P_n)$. So, we have proved
statement $3$.

To prove statement $1$, recall that the generators satisfy the
relations from the first statement. They are {\em defining}
relations since as it can be easily seen they guarantee that the
following equality holds,  $D=\sum Kx^\alpha \der^{[\beta]}  $,
where $D$ is an algebra generated by $x_i$, $\der_i^{[k]}$, $i=1,
\ldots , n$ and $k\geq 1$ that satisfy the given relations. Since
the algebra $\CD (P_n)$ is a factor algebra of $D$, the sum must
be
 a {\em direct} sum. Then the relations must be defining
 relations.

 Statement $4$ follows from statement $1$. $\Box $

\begin{theorem}\label{pCDAd=gd}
Let   $\i =(i_1, \ldots , i_r) \in \II_r$ and $\j =(j_1, \ldots ,
j_r)\in \JJ_r$, i.e. $\D =\D (\i , \j )\neq 0$, and  $\{ j_{r+1},
\ldots j_n\}=\{ 1, \ldots , n\} \backslash \{ j_1, \ldots , j_r\}$
and let  $A_\D$ be the localization of the algebra $A$ at the
powers of the element $\D$. Then
\begin{enumerate}
\item the algebra $\CD (A_\D )$ of differential operators on
$A_\D$ is a simple algebra such that

$$ \CD (A_\D )=\bigoplus_{k_{r+1}, \ldots , k_n\geq 0}A_\D
\d_{\i ; \j, j_{r+1}}^{[k_{r+1}]},\ldots ,\d_{\i ; \j,
j_{n}}^{[k_n]}=\bigoplus_{k_{r+1}, \ldots , k_n\geq 0} \d_{\i ;
\j, j_{r+1}}^{[k_{r+1}]},\ldots ,\d_{\i ; \j, j_{n}}^{[k_n]}A_\D$$
where $(\d_{\i ; \j, j_{s}}^{[k]})_{k\geq 0}$ is the unique
iterative higher derivation that is attached to the derivation $\D
(\i , \j )^{-1}\der_{\i ; \j, j_{s}}$, and all the elements
$\d_{\i ; \j, j_{s}}^{[k]}$, $k\geq 0$, $s=r+1, \ldots , n$
commute. Therefore, the algebra $\CD (A_\D )$ is generated by $A$
and $HS_K(A_\D )$.
 \item $\Der_K(A_\D )=\bigoplus_{\nu =r+1}^nA_\D
\der_{\i ; \j, j_{\nu}}$. \item For each $l\geq 0$,
$$ \CD (A_\D )_l=\bigoplus_{k_{r+1}+ \cdots + k_n\leq l}A_\D
\d_{\i ; \j, j_{r+1}}^{[k_{r+1}]},\ldots ,\d_{\i ; \j,
j_{n}}^{[k_n]}=\bigoplus_{k_{r+1}+ \cdots + k_n\leq l} \d_{\i ;
\j, j_{r+1}}^{[k_{r+1}]},\ldots ,\d_{\i ; \j,
j_{n}}^{[k_n]}A_\D.$$
\end{enumerate}
\end{theorem}

{\it Proof}. The second statement follows from Theorem
\ref{p25July04}.

Without loss of generality we can assume that $\i =(1,2, \ldots ,
r)$ and $\j =(1,2, \ldots , r)$. Let
$$\der_{r+1}:=\D^{-1}\der_{\i ; \j, r+1},\ldots  ,\der_n:=\D^{-1}\der_{\i ; \j,
n}. $$ Then $\der_i(x_j)=\d_{ij}$ for all $i,j=r+1, \ldots n $.
 By the  second statement, the commutator of the derivations
  $$[\der_i,\der_j]=\sum_{k=r+1}^na_{ij}^k\der_k\in \Der_K(A_\D ),
   \;\; a_{ij}^k\in A_\D ,$$
annihilates the elements $x_{r+1}, \ldots , x_n$. Therefore, all
$a_{ij}^k=0$ since
$$a_{ij}^k=\sum_{l=r+1}^na_{ij}^l\der_l (x_k) =[\der_i,\der_j]
(x_k) =0; $$ and it follows that the derivations $\der_i$ {\em
commute}.


The functor of taking derivations commutes with localizations,
therefore $\Der_K(A_\D )=\oplus_{i=r+1}^n A_\D \der_i$ implies
$\Der_K(Q)=\oplus_{i=r+1}^n Q \der_i$.  By \cite{BourbakiAlgII},
Theorem 1, Ch. V, Sec. 9, the field $Q$ is a finite  separable
extension of its subfield $Q':=K(x_{r+1}, \ldots x_n)$, hence
$\Der_{Q'}(Q)=0$. In characteristic $p>0$, a $p$'th power of a
$K$-derivation is again a $K$-derivation, hence $\der_i^p\in
\Der_K(A_\D )\subseteq \Der_K(Q)$. Since $Q'\subseteq \ker
(\der_i^p)$ and $Q$ is algebraic and separable over $Q'$, we must
have $\der_i^p=0$. Recall that over a perfect field any field
extension is separable (\cite{Ma}, Theorem 26.3). Now, by
\cite{Ma}, Theorem 27.4, each derivation $\der_i\in \Der_K(Q)$ can
be extended to an iterative higher derivation $\der_i^*:=({\rm
id}_Q, \der_i:= \der_i^{[1]}, \der_i^{[2]}, \ldots )\in HS_K(Q)$,
and it is  unique by \cite{Ma}, Theorem
 27.2. Now, considering the derivation $\der_i$ as an element of
 $\Der_K(Q')$, by the same arguments the derivation $\der_i$ can
 be uniquely extended to an iterative higher derivation
$$({\rm id}_{Q'}, \der_i, \frac{\der_i^2}{2!}, \ldots ,
\frac{\der_i^k}{k!}, \ldots )\in HS_K(Q'),$$
 then this iterative higher derivation has a {\em unique}
 extension to an iterative higher derivation of $Q$ (by \cite{Ma},
 Theorem 27.2). By uniqueness, it must coincide with $\der_i^*$,
 i.e. $\der_i^{[k]}|_{Q'}=\frac{\der_i^k}{k!}$ for all $i$ and
 $k$. A direct calculation gives,
 $$ -(\ad \, x_i) (\der_i^{[k]})=[\der_i^{[k]},
 x_i]=[\frac{\der_i^k}{k!}, x_i]=\frac{\der_i^{k-1}}{(k-1)!}=
\der_i^{[k-1]}.$$

Clearly, the inner derivations  $\ad (x_{r+1}), \ldots , \ad
(x_n)$ of the algebra $E:={\rm End}_K(A_\D )$ commute. Applying
several times Lemma \ref{pdx=1} we obtain  the algebra
$$ N=N(\ad (x_{r+1}), \ldots , {\rm ad}
(x_n); E):= \bigcap_{i=r+1}^nN(\ad (x_i), E)=\bigoplus_{\alpha \in
\mathbb{N}^{n-r}}C\der^{[\alpha ]} =\bigoplus_{\alpha \in
\mathbb{N}^{n-r}}\der^{[\alpha ]}C$$ where $\alpha =(\alpha_{r+1},
\ldots, \alpha_n)$, $\der^{[\alpha
]}:=\der_{r+1}^{[\alpha_{r+1}]}\cdots \der_{n}^{[\alpha_{n}]}$,
and $C:=\cap_{i=r+1}^n\ker (\ad (x_i))$ is the  subalgebra of $E$.
So, any element $u$ of $N$ is uniquely written as a sum
$u=\sum_{\alpha \in \mathbb{N}^{n-r}}c_\alpha \der^{[\alpha ]}$, $
c_\alpha \in C$. The algebra $N=\cup_{i\geq 0}N_i$ has a natural
filtration by the total degree of the $\der_i$'th, that is
$N_i=\bigoplus_{|\alpha |\leq i}C\der^{[\alpha ]}$ where $|\alpha
| := \alpha_{r+1}+\cdots +\alpha_n$. Clearly, $\CD (A_\D
)\subseteq N$ and $\CD (A_\D )_i\subseteq N_i$ for each $i\geq 0$.
Let us prove, by induction on $i$, that
$$\CD (A_\D
)_i=D_i:=\sum_{|\alpha|\leq i}A_\D \der^{[\alpha]}, \;\; i\geq
0.$$ The case $i=0$ is true, $\CD (A_\D )=A_\D =D_0$. Suppose that
$i>0$, and by induction $\CD (A_\D )_{i-1}=D_{i-1}$. Take $u\in
\CD (A_\D )_i$. Since $\CD (A_\D )_i\subseteq N_i$, the element
$u$ can be written a sum $u=\sum_{|\alpha |\leq i}c_\alpha
\der^{[\alpha ]}$ for some elements $c_\alpha \in C$. For each
$j=r+1, \ldots , n$,
$$-\ad (x_j)(u)=\sum_{|\alpha |\leq i}c_\alpha \der^{[\alpha
-e_j]}\in \CD(A_\D )_{i-1}=D_{i-1},$$
 (where the set $e_{r+1}, \ldots, e_n$ is the obvious `free basis' for
 $ \mathbb{N}^{n-r}:=\mathbb{N}e_{r+1}\oplus  \cdots \oplus \mathbb{N}e_n$),
 therefore all $c_\alpha \in
A_\D$ with $\alpha \neq 0$. Since $c_0=u-\sum_{\alpha \neq 0, |
\alpha | \leq i} c_\alpha \der^{[\alpha ]}\in C\cap \CD (A_\D
)_i$, it follows from the claim below that $c_0\in A_\D $.
Therefore, $\CD (A_\D )=D_i$, and so $\CD (A_\D )=A_\D \langle
\der_{r+1}^{[k]}, \ldots ,\der_n^{[k]}\rangle_{k\geq 1}$.

$\noindent $

 {\it Claim}. $C\cap \CD (A_\D )_i=A_\D $ {\em for all
} $i\geq 0$.

$\noindent $

We use induction on $i$. The case  $i=0$ is trivial, $C\cap \CD
(A_\D )_0=C\cap A_\D =A_\D $. Note that  the algebra $C$ is an
$A_\D $-bimodule, and so it is invariant under the inner
derivation $\ad (a)$ for each element $a \in A_\D $. If the
intersection $I_i:= C\cap \CD (A_\D )_i\neq A_\D$ for some $i\geq
1$, then obviously $I_1\neq A_\D $, and since
$$ \CD (A_\D )_1=A_\D +\Der_K(A_\D )=A_\D +\sum_{i=r+1}^nA_\D
\der_i,$$
 one can choose  an element $u=a_0+\sum_{i=r+1}^na_i\der_i\in I_1\backslash
 A_\D$ for some elements $a_i\in A_\D$ such that $a_j\neq 0$ for
 some $j\geq r+1$. We have a contradiction: $ 0=[u,x_j]=a_j\neq
 0$, which proves the claim.

 Next, let us  prove by induction on $s=k+l$ that the elements $\der_i^{[k]}$,
$\der_j^{[l]}$, $k,l\geq 0$, $i,j=r+1, \ldots , n$  commute. The
case $s=0$ is obviously true as $\der_i^{[0]}={\rm id}_{A_\D}$ for
all $i$. Suppose that $s>0$ and the result is true for all $s'<s$.
For each $t=r+1, \ldots , n$,
$$ [[\der_i^{[k]},
\der_j^{[l]}], x_t]=[[\der_i^{[k]}, x_t], \der_j^{[l]}]
+[\der_i^{[k]},[ \der_j^{[l]}, x_t]]=\d_{it}[\der_i^{[k-1]},
\der_j^{[l]}]+\d_{jt}[\der_i^{[k]}, \der_j^{[l-1]}]=0,$$ by
induction. By the claim, $[\der_i^{[k]}, \der_j^{[l]}]\in A_\D$,
and so $[\der_i^{[k]}, \der_j^{[l]}]=[\der_i^{[k]},
\der_j^{[l]}]*1=0$.

To prove that the algebra $\CD (A_\D )$ is simple, let $L$ be a
nonzero ideal of $\CD (A_\D )$. It remains to prove that $L=\CD
(A_\D )$. Take a nonzero element, say $u$, of $L$. Applying
several times maps of the type $\ad (x_{j_k})$, $r+1\leq k\leq n$,
to the element $u$ we have a nonzero element, say $u_1\in L\cap
A_\D$. Since $\Kdim (A_\D )=\Kdim (K[x_{r+1}, \ldots , x_n])$, we
must have $A_\D u_1\cap K[x_{r+1}, \ldots , x_n]\neq 0$. Pick a
nonzero element, say $u_2$, of the intersection, then $u_2\in \CD
(K[x_{r+1}, \ldots , x_n])\subseteq \CD (A_\D )$. The algebra $\CD
(K[x_{r+1}, \ldots , x_n])$ is a simple algebra (Corollary
\ref{c1pdx=1}), hence $L=\CD (A_\D )$. This proves that the
algebra $\CD (A_\D )$ is a simple algebra. The rest is obvious.
 $\Box $

It is clear from Theorem \ref{pCDAd=gd}, that  a differential
operator from the algebra $\CD (A_\D )$ (or from the algebra $\CD
(A)\subseteq \CD (A_\D )$) is uniquely determined by its action on
the polynomial subalgebra $P:= K[x_{j_{r+1}}, \ldots,  x_{j_n}]$,
that is the restriction map 
\begin{equation}\label{DresP}
\CD (A_\D )\ra \Hom_K(P, A_\D ), \;\; \d \mapsto \d|_P,
\end{equation}
 is an injective map. Then the second statement of Corollary
 \ref{1pCDAd=gd} follows, the first statement of Corollary
 \ref{1pCDAd=gd} was already proved in the proof of Theorem
 \ref{pCDAd=gd}.

\begin{corollary}\label{1pCDAd=gd}
Let  $\i = (i_1, \ldots , i_r)\in \II_r$ and $\j = (j_1, \ldots ,
j_r)\in \JJ_r$, i.e. $\D (\i , \j )\neq 0$, and $\{ j_{r+1},
\ldots , j_n\} = \{ 1, \ldots , n \} \backslash \{ j_1, \ldots ,
j_r\}$. Then
\begin{enumerate}
\item The field $Q$ is a finite separable field extension of the
field $K(x_{j_{r+1}}, \ldots ,x_{j_n})$. \item The derivations
$\D(\i , \j )^{-1}\der_{\i ; \j , j_{r+1}}, \ldots , \D(\i , \j
)^{-1}\der_{\i ; \j , j_n}$ from Theorem \ref{pCDAd=gd} are
respectively the partial derivatives
$\der_{j_{r+1}}:=\frac{\der}{\der x_{j_{r+1}}}, \ldots ,
\der_{j_n}:=\frac{\der}{\der x_n}$ of the algebra $A_\D $ (and of
the field of fractions of $A$).
\end{enumerate}
\end{corollary}

{\it Remarks}. 1. Statement 1 of Corollary
 \ref{1pCDAd=gd} is a strengthening of the following well known
 result -- Theorem 26.2,  \cite{Ma} (and the Remark after Theorem 26.2,  \cite{Ma}): Let $K$ be a field of characteristic
 $p$ and $L= K(y_1, \ldots , y_t)$ be a finitely generated field
 extension of $K$ which is separably generated over $K$; then
 there exists a subset of the set of generators, say $y_{i_1},
 \ldots ,y_{i_d}$, $d=\trdeg_K(L)$, such that $L$ is separable over
 $K(y_{i_1}, \ldots ,y_{i_d})$.

2. The equality $\der_{j_{r+k}}:=\frac{\der}{\der x_{j_{r+k}}}$
means that the derivation $\der_{j_{r+1}}$ is a unique extension
of the partial derivative $\frac{\der}{\der x_{j_{r+1}}}$ of the
polynomial algebra  $K[x_{r+1}, \ldots , x_n]$ to the algebra
$A_\D$.

\begin{theorem}\label{rpCDAd=gd}
Keep the assumption of Theorem \ref{pCDAd=gd}. Then the algebra
$\CD (A_\D )$ is generated by the algebra $A_\D$ and the elements
 $\d_{\i ; \j, j_{\nu}}^{[k]}$ where   $\nu =r+1, \ldots , n$ and $k\geq 0$,
that satisfy the following defining relations: for all $\nu , \mu
=r+1, \ldots , n$ and $k,l\geq 1$ (where $\d^{[0]}_{\i , \j ,
j_\nu }:=1$)
$$[\d_{\i ; \j, j_\nu }^{[k]}, \d_{\i ; \j, j_\mu }^{[l]}]=0,
\;\; \; \;\;\;  \d_{\i ; \j, j_\nu }^{[k]}\d_{\i ; \j, j_\nu
}^{[l]}={k+l\choose k} \d_{\i ; \j, j_\nu }^{[k+l]},
\;\;\;\;\;\;\;  [\d_{\i ; \j, j_\nu }^{[k]}, x_{j_{\mu}}
]=\d_{j_\nu ,j_\mu} \d_{\i ; \j, j_\nu }^{[k-1]},
$$
and, for all $\nu =r+1, \ldots , n$ and $s=1, \ldots , r$,
$$ [\d_{\i ; \j, j_\nu }^{[k]}, x_{j_s}
]= \sum_{t=1}^k  \d_{\i ; \j, j_\nu }^{[t]}(x_{j_s})\d_{\i ; \j,
j_\nu }^{[k-t]},
$$
where  $\d_{\i ; \j, j_\nu }^{[t]}(x_{j_s})= \frac{(\D^{-1}
\der_{\i , \j , j_\nu})^t}{t!}(x_{j_s})\in A_\D $.
\end{theorem}

{\it Remark}. The elements $\d_{\i ; \j, j_{\nu}}^{[t]}(x_{j_s})$
can be found explicitly by combining Corollary \ref{1pCDAd=gd}.(1)
and Theorem \ref{12Jul8}.

{\it Proof}. Clearly, the generators satisfy the  given relations.
Suppose that $D$ is an algebra generated by the given elements
that satisfy the given defining relations. One can easily see that
$$ D=\sum_{k_{r+1}, \ldots , k_n\geq 0} A_\D
\d_{\i ; \j, j_{r+1}}^{[k_{r+1}]},\ldots ,\d_{\i ; \j,
j_{n}}^{[k_n]}.$$ By Theorem \ref{pCDAd=gd}.(1), the sum above
must be direct. Therefore, $D=\CD (A_\D )$, which implies that the
relations are defining relations for the algebra $\CD (A_\D )$.
$\Box $

\begin{corollary}\label{2pCDAd=gd}
 Let $Q$ be the field of fractions of the algebra $A$. Under the
 assumption of Theorem \ref{pCDAd=gd},
\begin{enumerate}
\item the algebra $\CD (Q )$ of differential operators on $Q$ is a
simple algebra such that

$$ \CD (Q)=\bigoplus_{k_{r+1}, \ldots , k_n\geq 0}Q
\d_{\i ; \j, j_{r+1}}^{[k_{r+1}]},\ldots ,\d_{\i ; \j,
j_{n}}^{[k_n]}=\bigoplus_{k_{r+1}, \ldots , k_n\geq 0} \d_{\i ;
\j, j_{r+1}}^{[k_{r+1}]},\ldots ,\d_{\i ; \j, j_{n}}^{[k_n]}Q$$
where $(\d_{\i ; \j, j_{s}}^{[k]})_{k\geq 0}\in HS_K(Q )$ is the
iterative higher derivation which is a unique extension of the
$(\d_{\i ; \j, j_{s}}^{[k]})_{k\geq 0}\in HS_K(A_\D )$ from
Theorem \ref{pCDAd=gd}, they also commute. Therefore, the algebra
$\CD (Q)$ is generated by the field $Q$ and $\HS_K(Q)$.
 \item $\Der_K(Q)=\bigoplus_{\nu =r+1}^nQ
\der_{\i ; \j, j_{\nu}}$. \item For each $l\geq 0$,
$$ \CD (Q )_l=\bigoplus_{k_{r+1}+ \cdots + k_n\leq l}Q
\d_{\i ; \j, j_{r+1}}^{[k_{r+1}]},\ldots ,\d_{\i ; \j,
j_{n}}^{[k_n]}=\bigoplus_{k_{r+1}+ \cdots + k_n\leq l} \d_{\i ;
\j, j_{r+1}}^{[k_{r+1}]},\ldots ,\d_{\i ; \j, j_{n}}^{[k_n]}Q.$$
\end{enumerate}
\end{corollary}

{\it Proof}. Note that $\CD (Q)\simeq Q\t_{A_\D}\CD(A_\D )$ and
$\Der_K(Q)\simeq Q\t_{A_\D }\Der_K(A_\D )$, and the results follow
from Theorem \ref{pCDAd=gd}.  $\Box $

\begin{theorem}\label{3Feb05}
Let the algebra $A$ be a regular algebra. Then the algebra $\CD
(A)$ of differential operators on $A$ is a simple algebra
generated by $A$ and $HS_K(A)$, that is $\CD (A)=\D (A)$.
\end{theorem}

{\it Proof}. Let $\D =\D (A)$ be the subalgebra of ${\rm
End}_K(A)$ generated by $A$ and $HS_K(A)$. By Theorem
\ref{pCDAd=gd}, $\CD (A)_{\D (\i , \j )}=\D_{\D (\i , \j )}$ for
all non-singular $\i$ and $\j $, or equivalently $\prod_{\i , \j
}\CD (A)_{\D (\i , \j )}\t_{\CD (A)}(\CD (A)/\D ) =0$. By Theorem
\ref{CDAffl}, we must have $\CD (A)=\D $. By Theorems \ref{CDAffl}
and \ref{pCDAd=gd}, $\CD (A)$ is a simple algebra. $\Box $

\begin{proposition}\label{25Jul8}
Let the algebra $A$ be a regular algebra. Then $\CD (A)=
\bigcap_{\i \in \II_r, \j\in\JJ_r}\CD (A)_{\D (\i, \j )}$ where
the intersection is taken in the algebra $\CD (Q)$.
\end{proposition}

{\it Proof}. We denote by $\CD'$ the intersection. Then $\CD
(A)\subseteq \CD'$ and $ \CD (A)_{\D (\i, \j )}= \CD'_{\D (\i, \j
)}$ for all $\i \in \II_r$ and $ \j\in\JJ_r$. Then $\CD (A) =\CD'$
since the extension $\CD (A) \ra \prod_{\i \in \II_r,
\j\in\JJ_r}\CD (A)_{\D (\i, \j )}$ is faithfully flat. $\Box$.

\begin{theorem}\label{12Jul8}
Let $k\subseteq L\subseteq L'$ be fields such that $L'= L(x)$ for
a separable element $x\in L'$ over $L$ and $f(t) =
t^s+a_{s-1}t^{s-1}+\cdots + a_0\in L[t]$ be a minimal polynomial
for $x$. Then each higher derivation $\{ \d^{[i]}, i\geq 0\}\in
\HS_k(L)$ can be uniquely extended to a higher derivation $\{
\d^{[i]}, i\geq 0\}\in \HS_k(L')$. Moreover, 
\begin{equation}\label{extdix}
\d^{[i]}(x) = -\frac{1}{f'(x)}\sum_{m=0}^s\; \sum_{i_0+\cdots +
i_m=i, i_\nu \neq i , \nu \geq 1}
\d^{[i_0]}(a_m)\d^{[i_1]}(x)\cdots \d^{[i_m]}(x), \;\; i\geq 1,
\end{equation}
where $f'=\frac{df}{dx}$.
\end{theorem}

{\it Proof}. The $k$-algebra homomorphism $\s (\cdot ) =
\sum_{i\geq 0} \d_i(\cdot ) t^i: L\ra L[[t]]$ can be extended to a
$k$-algebra homomorphism $\s (\cdot ) = \sum_{i\geq 0} \d_i(\cdot
) t^i: L'\ra L'[[t]]$ iff $0= \s (f(x)) = \sum_{i\geq 0}
\d_i(f(x)) t^i$ iff $\d_i(f(x))=0$, $i\geq 1$ iff the equalities
of the theorem hold (for each $i\geq 1$, the equality
$\d_i(f(x))=0$ can be rewritten as (\ref{extdix})). $\Box $

\begin{corollary}\label{a12Jul8}
Let $k\subseteq L\subseteq L'$ be fields such that $L'$ is
separable over $L$. Then $\HS_k(L)\subseteq \HS_k(L')$.
\end{corollary}

\begin{corollary}\label{b12Jul8}
Let $K\subseteq Q\subseteq Q'$ be fields finitely generated over
the field $K$ and let  the field $Q'$ be separable over $Q$. Then
$\CD_K(Q') = Q'\CD_K(Q)$.
\end{corollary}

{\it Proof}. We may assume that the field $Q$ is as in Corollary
\ref{2pCDAd=gd} (by inverting, if necessary, one of the nonzero
minors $\D (\i , \j )$). By Corollary \ref{1pCDAd=gd}.(1), the
field $Q$ is a separable extension of its subfield $L:=
K(x_{j_{r+1}}, \ldots , x_{j_n})$, then so is $Q'$ since $Q'$ is
separable over $Q$. The higher derivation $\{ \d_{\i ; \j ,
j_s}^{[k]}, k\geq 0\}\in \HS_K(Q)$ can be extended uniquely to a
higher derivation of the field $Q'$, by Corollary \ref{a12Jul8}.
When we write down statement 1 of Corollary \ref{2pCDAd=gd} for
the fields $Q$ and $Q'$ the equality $\CD_K(Q') = Q'\CD_K(Q)$
follows at once.  $\Box $

$\noindent $

Note that in Corollary \ref{b12Jul8}  the inclusion
$\CD_k(Q)\subseteq \CD_K(Q')$ follows from Corollary
\ref{2pCDAd=gd}.(1), Corollary \ref{a12Jul8} and the simplicity of
the algebra $\CD_k(Q)$. Then the inclusion $Q'\CD_K(Q)\subseteq
\CD_K(Q')$ is obvious.

\begin{proposition}\label{P3Feb05}
Let $\i , \i'\in \II_r$, $\j =(j_1, \ldots , j_r)\in \JJ_r$, $\j'
=(j_1', \ldots , j_{r+1}')\in \JJ_{r+1}$, and $\{ j_{r+1}, \ldots
j_n\} =\{ 1, \ldots , n\} \backslash \{ j_1, \ldots , j_r\}$. Then
\begin{enumerate}
 \item

\begin{equation}
 \der_{\i',\j'}=\D (\i , \j )^{-1}\sum_{l=1}^s(-1)^{r+1+\nu_l}\D (\i';
j_1', \ldots , \widehat{j_{\nu_l}'}, \ldots , j_{r+1}')\der_{\i
;\j , j_{\nu_l}'}
\end{equation}
where  $j_{\nu_1}', \ldots , j_{\nu_s}'$ are  the elements
 of the set $\{ j_1', \ldots , j_r', j_{r+1}'\} \backslash \{ j_1,
 \ldots , j_r\}$.
\item $\der_{\i',\j'}=(-1)^{r+1+k}\frac{\D (\i'; j_1', \ldots ,
\widehat{j_k'}, \ldots , j_{r+1}' )}{\D (\i ; j_1', \ldots ,
\widehat{j_k'}, \ldots , j_{r+1}'  )}\der_{\i ,\j'}$ provided $\D
(\i ; j_1', \ldots , \widehat{j_k'}, \ldots , j_{r+1}'  )\neq 0$.
\end{enumerate}
\end{proposition}

{\it Proof}. $1$. By Theorem \ref{p25July04},  $\der_{\i' ,
\j'}=\D (\i , \j )^{-1}\sum_{k=r+1}^n\l_k\der_{\i ; \j , j_k}$ for
some $\l_k\in A$. For each $k'=r+1, \ldots , n$, $\der_{\i ; \j ,
j_k}(x_{j_{k'}})=\d_{k,k'}\D (\i , \j )$. Then evaluating the
equality above at $x_{j_k}$, we
 get the equality $\l_k= \der_{\i' ,\j'}(x_{j_k})$. So, $\l_k=0$
 if $j_k\not\in \{ j_{\nu_1}', \ldots , j_{\nu_s}'\}$, and if
 $j_k=j_{\nu_l}'$ then $\l_k=(-1)^{r+1+\nu_l}
 \D (\i ; j_1', \ldots , \widehat{j_{\nu_l}'}, \ldots ,
 j_{r+1}')$. This finishes the proof of the first statement.

 $2$. By the first statement where we put $\j =
 (j_1', \ldots , \widehat{j_k'}, \ldots ,
 j_{r+1}')$, we have $\der_{\i',\j'}=(-1)^{r+1+k}\frac{\D (\i'; j_1', \ldots ,
\widehat{j_k'}, \ldots , j_{r+1}' )}{\D (\i ; j_1', \ldots ,
\widehat{j_k'}, \ldots , j_{r+1}'  )}\der_{\i ,\j'}$.
 $\Box $

$\noindent $

{\it Remark}. Let us fix elements $\i\in \II_r$ and $\j\in \JJ_r$.
Then, for each $\i'\in \II_r$ and $\j'\in \JJ_{r+1}$ (as above),
let $a(\i',\j')$ be the vector of coefficients $(\l_{r+1}, \ldots
, \l_n)$ from the proof of Proposition \ref{P3Feb05}. If $A$ is a
{\em regular} algebra, the vectors $a(\i',\j')$ form a generating
set for the $A$-module of solutions to the system of inclusions
from Theorem \ref{p25July04} (by Theorem \ref{9bFeb05}).

$\noindent $

{\bf  Proof of Theorem \ref{21Jul8}}. Let $D=D(A)$ be the algebra
generated by the algebra $A$ and the elements $d_{\i ;  \j ,
j_\nu}^{[k]}$ that satisfy the defining relations (R1)--(R5).
Theorem \ref{21Jul8} follows easily from the Claim, and the Claim
follows from Theorem \ref{pCDAd=gd} and Theorem \ref{rpCDAd=gd}.

$\noindent $

{\it Claim: There is an algebra homomorphism $D\ra \CD := \CD (A)$
such that $D_{\D (\i , \j )}\simeq \CD_{\D (\i , \j )}$ for all
$\i \in \II_r$ and $\j \in \JJ_r$ where $D_{\D (\i , \j )}$ is the
(left and right) Ore localization of the algebra $D$ at the powers
of the element $\D (\i , \j )$. }

$\noindent $

Indeed, we have the commutative diagram of algebra homomorphisms
$$\xymatrix{D\ar[r]\ar[d]&\prod_{\i\in \II_r , \j \in \JJ_r}D_{\D (\i , \j )}\ar[d]\\
\CD \ar[r]&\prod_{\i \in \II_r, \j\in \JJ_r }\CD_{\D (\i , \j )}}
$$
where the horizontal maps are faithfully flat extensions (since
$A=(\D (\i, \j ))_{\i\in \II_r, \j \in \JJ_r}$ as $A$ is regular)
and the right vertical map is an isomorphism, by the Claim. By
faithful flatness, the left vertical map is an isomorphism, i.e.
$D\simeq \CD$, and the ring $\CD$ is simple since each ring
$\CD_{\D (\i , \j )}$ is so.

It remains to prove the Claim. By Theorem \ref{rpCDAd=gd}, for the
elements $\i \in \II_r$ and $\j \in \JJ_r$, the algebra $\CD_\D$
(where $\D := \D (\i , \j )$) is generated by the algebra $A_\D$
and the elements $\{ \d_{\i ; \j , j_\nu}^{[k]}, k\geq 0\}$, $\nu
= r+1, \ldots , n$, that satisfy the four types of the defining
relations of Theorem \ref{rpCDAd=gd}. It is obvious that the
algebra $\CD_\D$ is generated by the algebra $A_\D$ and the
elements $\{ d_{\i ; \j , j_\nu}^{[k]}, k\geq 0 \}$, $\nu = r+1,
\ldots , n$. When we multiply the fourth relation of Theorem
\ref{rpCDAd=gd} by the invertible element $\D^{n(k) +k}$ on the
left we obtain the relation (R4) which is equivalent to the
original one. Using the relation (R4), we see that
\begin{equation}\label{R4Ds}
\D^{-s}\d_{\i ; \j , j_\nu}^{[k]}\D^{n(k)+k+s}\in \sum_{t=0}^kA
d_{\i ; \j , j_\nu}^{[t]}, \;\; s\geq 0.
\end{equation}
In more detail, let $d^{[t]}:=d_{\i ; \j , j_\nu}^{[t]}$. By (R4),
$(\ad (\D ))^j (d^{[k]})\in \sum_{l=0}^{k-j}Ad^{[l]}$. Then
\begin{eqnarray*}
 \D^{-s}\d_{\i ; \j , j_\nu}^{[k]}\D^{n(k)+k+s}&=& \D^{-s-n(k)}d^{[k]}\D^{n(k)+k+s}\\
 &=& \D^{-s-n(k)} \sum_{j=0}^k {n(k)+k+s\choose j} \D^{n(k)
 +k+s-j} (-\ad (\D ))^j (d^{[k]})\\
 &\in &\sum_{t=0}^k Ad^{[k]}.\\
\end{eqnarray*}
By (\ref{R4Ds}), 
\begin{equation}\label{R4Ds1}
\D^{-s}\d_{\i ; \j , j_\nu}^{[k]}\D^{n(k)+k+s}= \sum_{t=0}^k a_{\i
; \j , j_\nu}^{[k]}(s,t)d_{\i ; \j , j_\nu}^{[t]}, \;\; s\geq 0,
\end{equation}
for some elements $a_{\i ; \j , j_\nu}^{[k]}(s,t)\in A$. By
Theorem 3.1.(5), \cite{JPAA-08invcharp}: 
\begin{equation}\label{R4Ds2}
a_{\i ; \j , j_\nu}^{[k]}(s,t)= \psi_{j_\nu}((-\ad\, x_{j_\nu})^t
(\D^{-s} \d_{\i ; \j , j_\nu}^{[k]}\D^{n(k)+k+s}))\cdot
\D^{-n(t)},
\end{equation}
where $\psi_{j_\nu}(\cdot ) := \sum_{j\geq 0} (\ad \,
x_{j_\nu})^j(\cdot ) \d_{\i ; \j , j_\nu }^{[j]}:\CD_\D \ra
\CD_\D$. When we multiply the first relation of Theorem
\ref{rpCDAd=gd} by the element $\D^{n(k)+n(l)}$ on the left and by
the element $\D^{n(k)+n(l)+k+l}$ on the right we obtain the
relation (R1) using (\ref{R4Ds1}). Similarly, multiplying the
second relation of Theorem \ref{rpCDAd=gd} by the element
$\D^{n(k+l)}$ on the left and by the element $\D^{n(l)+l}$ on the
right we obtain the relation (R2) using (\ref{R4Ds1}). When we
multiply the third relation of Theorem \ref{rpCDAd=gd} by the
element $\D^{n(k)}$ on the left we obtain the relation (R3).

It is obvious that the algebra $\CD_\D$ is generated by the
algebra $A_\D$ and the elements $\{ d_{\i ; \j , j_\nu}^{[k]},
k\geq 0\}$, $\nu = r+1, \ldots , n$ that satisfy the defining
relations (R1)--(R4). Recall that, for  $\i'\in \II_r$ and $\j'\in
\JJ_r$, we have $\D':= \D(\i', \j')$. By Theorem \ref{pCDAd=gd},
\begin{equation}\label{blk1}
\d_{\i'; \j', j'_\s}^{[l]}= \sum b_{l, \k } \d_{\i; \j,
j_{r+1}}^{[k_{r+1}]}\cdots \d_{\i; \j, j_n}^{[k_n]}, \;\; l\geq 1,
\end{equation}
for some elements $b_{l, \k }=b_{l, \k }(\i'; \j', j'_\s )\in
A_\D$ where the sum is taken over the vectors $\k := (k_{r+1},
\ldots , k_n)\in \N^{n-r}$ such that $|\k | \leq l$. A formula for
the elements $b_{l, \k }$ is given by Theorem 3.1.(5),
\cite{JPAA-08invcharp}: 
\begin{equation}\label{blk2}
b_{l,\k }= \psi_{j_{r+1}}\cdots \psi_{j_n}\prod_{\nu = r+1}^n
(-\ad \, x_{j_\nu })^{k_{r+\nu}} (\d^{[l]}_{\i'; \j', j'_\s})
\end{equation}
where $\psi_{j_\nu}(\cdot ) := \sum_{j\geq 0} (\ad \,
x_{j_\nu})^j(\cdot ) \d_{\i ; \j , j_\nu }^{[j]}:\CD_\D \ra
\CD_\D$.

For each natural number $l$, choose a natural number $m(l)$ such
that 
\begin{equation}\label{clk1}
c_{l, \k } := \D^{m(l)}b_{l,\k }\in A
\end{equation}
for all $\k \in \N^{n-r}$ with $|\k | \leq l$, and all  $\i ,
\i'\in \II_r$, $\j , \j'\in \JJ_r$ and $j'_\s$. When we multiply
the equality (\ref{blk1}) by the element $\D^{m(l)}\D'^{n(l)}$ on
the left and by the element $\D^{\S_{n-r}}$ on the right and make
some transformations using (\ref{R4Ds1}) several times we obtain
the relation (R5). In more detail, let $\d^{[k_s]}:=
\d^{[k_s]}_{\i ; \j , j_\nu}$. Then
\begin{eqnarray*}
 \D^{m(l)}d_{\i'; \j', j'_\s}^{[l]}\D^{\S_{n-r}}&=& \D'^{n(l)}\sum c_{l,\k} \d^{[k_{r+1}]}\cdots \d^{[k_n]} \cdot \D^{\S_{n-r}}\\
 &=&\D'^{n(l)}\sum c_{l,\k}\D^{-\S_0} \d^{[k_{r+1}]}\D^{\S_1}\cdot \D^{-\S_1}\d^{[k_{r+2}]}\D^{\S_2}
 \cdots  \D^{-\S_{n-r-1}}\d^{[k_n]}\D^{-\S_{n-r}}\\
 & \cdot& \D^{-\S_{n-r}+\S_{n-r}}=\D'^{n(l)}\sum
 c_{l,\k}\prod_{\rho =1}^{n-r} \D^{-\S_{\rho
 -1}}\d^{[k_{r+\rho}]}\D^{\S_\rho}\\
 &=& \D'^{n(l)}\sum
 c_{l,\k}\prod_{\rho =1}^{n-r} \sum_{t_{r+\rho}=0}^{k_{r+\rho}}a_{\i ; \j
 , j_\nu}^{[k_{r+\rho}]}(\S_{\rho -1}, t_{r+\rho}) d_{\i ; \j
 , j_\nu}^{[t_{r+\rho}]}.
\end{eqnarray*}
The relations (R3) guarantee that the (left and right) Ore
localization $D_\D$ exists. Now, the Claim is obvious where $D\ra
\CD$, $d_{\i ; \j , j_\nu}^{[k]}\mapsto d_{\i ; \j
 , j_\nu}^{[k]}$ is the natural algebra homomorphism over $A$. The
 proof of Theorem \ref{21Jul8} is complete. $\Box$

$\noindent $

{\bf  Proof of Theorem \ref{9bFeb05}}.  Let ${\rm DER}(A)$ be a
left  $A$-module generated by symbols $\der_{\i , \j }$ subject to
the defining relations (\ref{Derel}). We have the  commutative
diagram of left $A$-modules:
$$\xymatrix{{\rm DER}(A)\ar[r]\ar[d]&\prod_{\i\in \II_r , \j \in \JJ_r}{\rm DER}(A)_{\D (\i , \j )}\ar[d]\\
\Der_K(A) \ar[r]&\prod_{\i \in \II_r, \j\in \JJ_r }\Der_K(A)_{\D
(\i , \j )}}
$$
where the horizontal maps are faithfully flat $A$-module
monomorphisms
 as  $A=(\D (\i , \j ))_{\i\in \II_r, \j \in \JJ_r}$
 ($A$ is regular) and the vertical maps are natural $A$-module
 epimorphisms. By Proposition  \ref{P3Feb05}, Theorem \ref{p25July04},  and (\ref{Derel}),
 each epimorphism ${\rm DER}(A)_{\D (\i , \j )}\ra \Der_K(A)_{\D (\i ,
 \j )}$ is an isomorphism. So, the right vertical map must be an
 isomorphism, and so the left vertical map must be an
 isomorphism (by faithfully flatness), i.e. ${\rm DER}(A)\simeq \Der_K(A)$.   $\Box $

\begin{theorem}\label{24Jul8}
The set $\HS_K(A)$ of higher derivations of the algebra $A$ leaves
invariant the Jacobian ideal $\ga_r$ of the algebra $A$.
\end{theorem}

{\bf Proof of Theorem \ref{23Jul8}}. $(1\Rightarrow 2)$  Theorem
\ref{3Feb05}.

$(2\Rightarrow 3 )$ Suppose that the $\D (A)$-module $A$ is not
simple then it contains a proper ideal, say $\ga$, stable under
$\HS_K(A)$. Then $\ga \D (A)$ is a proper  ideal of the algebra
$\D (A)$ since $0\neq \ga \D (A)(A)\subseteq \ga$, a
contradiction.

$(3\Rightarrow 1)$ By Theorem  \ref{24Jul8}, the Jacobian ideal
$\ga_r$ is a nonzero $\D (A)$-submodule of $A$, therefore
$A=\ga_r$ since $A$ is a simple $\D (A)$-module. So, $A$ is a
regular algebra. $\Box $


\section{Ring of differential operators
 on singular irreducible affine algebraic
variety}\label{opsing}

In this Section, we prove Theorem \ref{pCDAifg} (see Proposition
\ref{psinDAfg}.(2)), the local finiteness of the ring $\CD (A)$ of
differential operators  on the algebra $A$.

 \begin{lemma}\label{duchi}
Let $R=K\langle x_1, \ldots , x_n\rangle$ be a commutative
finitely generated algebra over the field $K$, and $\CD (R)$ be
the ring of differential operators on $R$. Each element $\d \in
\CD (R)_i$ is completely determined by its values on the elements
$x^\alpha$, $\alpha \in \mathbb{N}^n$, $| \alpha |\leq i$.
\end{lemma}

{\it Proof}. It suffices to prove that if an element  $\d \in \CD
(R)_i$ satisfies  $\d (x^\alpha )=0$ for all $\alpha$ such that $|
\alpha | \leq i$ then $\d =0$. We use induction on $i$. The case
$i=0$ is trivial: $\d \in \CD (R)_0=R$ and $0=\d \cdot 1=\d$.
Suppose that $i\geq 1$ and the statement is true for all $i'<i$.
For each $x_j$, $[\d , x_j]\in \CD (R)_{i-1}$ and, for each
$x^\alpha$ with $|\alpha | \leq i-1$, $[\d , x_j](x^\alpha )=\d
(x_jx^\alpha)-x_j\d (x^\alpha)=0$. By induction, $[\d , x_j]=0$.
Now, for any $x^\alpha$,
$$ \d (x^\alpha )=\d (x_jx^{\alpha -e_j})=x_j\d (x^{\alpha
-e_j})+[\d , x_j](x^{\alpha -e_j})=x_j\d (x^{\alpha -e_j})=\cdots
= x^\alpha \d (1)=0,$$
 and so $\d =0$, as required.
$\Box $

Let $S_1$ be a multiplicatively closed subset of the algebra $A$.
Let us consider a natural inclusion  $\CD (A)\subseteq S_1^{-1}
\CD (A)$ of filtered algebras (by the total degree of
derivations). $\CD (A)=\{ \d \in S_1^{-1} \CD (A)\, | \, \d
(A)\subseteq A\}$ and $\CD (A)_i=\{ \d \in S_1^{-1} \CD (A)_i\, |
\, \d (A)\subseteq A\}$, $i\geq 0$.

\begin{lemma}\label{1duchi}
Let  $\d \in S_1^{-1} \CD (A)_i$. Then $\d \in \CD (A)$ iff $\d
(x^\alpha )\in A$ for all $\alpha$ such that $| \alpha |\leq i$.
\end{lemma}

{\it Proof}. $(\Rightarrow )$ Trivial.

$(\Leftarrow )$  We use induction on $i$. When  $i=0$, $\d \in
S_1^{-1} \CD (A)_0=S_1^{-1} A$ and $\d =\d (1)=\d \cdot 1\in A$.
Suppose that $i\geq 1$ and the statement is true for all $i'<i$.
Let $\d \in S_1^{-1} \CD (A)_i$ satisfy $\d (x^\alpha )\in A$ for
all $\alpha $ such that $| \alpha | \leq i$. For each $j$,
 $[\d ,
x_j]\in S_1^{-1} \CD (R)_{i-1}$ and, for each $x^\alpha$ with
$|\alpha | \leq i-1$, $[\d , x_j](x^\alpha )=\d
(x_jx^\alpha)-x_j\d (x^\alpha)\in A$, and so, by induction, $[\d ,
x_j]\in \CD (A)_{i-1}$. Now, for any $x^\alpha$,
\begin{eqnarray*}
\d (x^\alpha )&=&\d (x_jx^{\alpha -e_j})=x_j\d (x^{\alpha
-e_j})+[\d , x_j](x^{\alpha -e_j})\equiv x_j\d (x^{\alpha -e_j})\,
{\rm mod} \, A\\
&\equiv &\cdots \equiv  x^\alpha \d (1)\equiv 0 \, {\rm mod} \, A,
\end{eqnarray*}
 and so $\d \in \CD (A)_i$.
$\Box $

\begin{proposition}\label{psinDAfg}
 Let  $\i =(i_1, \ldots , i_r)\in \II_r$, $\j =(j_1, \ldots , j_r)\in \JJ_r$, $\{ 1, \ldots , n\}
\backslash \{ j_1, \ldots , j_r\}=$ $ \{ j_{r+1}, \ldots , j_n\}
$, and we keep the notation of Theorem \ref{pCDAd=gd}. Then
\begin{enumerate}
\item For each $i\geq 0$, $\CD (A)_i=\{ \d \in \sum_{| \k  |\leq
i} A\d^{[\k ]} \, | \,  \d (x^\beta )\in A$ for all $\beta \in
\mathbb{N}^n$ such that $|\beta | \leq i\}$  where $\k =(k_{r+1},
\ldots , k_n)\in \mathbb{N}^{n-r}$ and  $\d^{[\k ]}
:=\prod_{s=r+1}^n\d_{\i ; \j, j_s}^{[k_s]}$ . \item For each
$i\geq 0$, $\CD (A)_i$ is a finitely generated left $A$-module.
\end{enumerate}
\end{proposition}

{\it Proof}. $1$. Clearly, $\CD (A)_i=\{ \d \in \sum_{|\k |\leq
i}A_\D \d^{[\k ]}\, | \, \d (A)\subseteq A\}$. Then, by Lemma
\ref{1duchi}, $\CD (A)_i=\{ \d \in \sum_{|\k |\leq i}A_\D \d^{[\k
]}\, | \, \d (x^\beta)\subseteq A$ for all $\beta \in
\mathbb{N}^n$ such that $|\beta | \leq i\}$. Then the conditions
that $\d (x_{j_{r+1}}^{\gamma_{r+1}}\cdots x_{j_n}^{\gamma_n})\in
A$ for all $\gamma=(\gamma_{r+1}, \ldots ,\gamma_n)\in
\mathbb{N}^{n-r}$ with $|\gamma | \leq i$ are equivalent to $\d
\in \sum_{|\k | \leq i}A\d^{[\k ]}$. This gives the first
statement.

$2$. $\CD (A)_i$ is the Noetherian left $A$-module as a submodule
of the Noetherian $A$-module $\sum_{|\k | \leq i}A\d^{[\k ]}$, and
so $\CD (A)_i$ is a finitely generated left $A$-module. $\Box $

Department of Pure Mathematics

University of Sheffield

Hicks Building

Sheffield S3 7RH

UK

email: v.bavula@sheffield.ac.uk


\begin{thebibliography}{99}

\bibitem{Alv-Blickle-Lyub} J. Alvarez-Montaner, M. Blickle, G. Lyubeznik,
 Generators of D-modules in positive characteristic,  {\em Math. Res. Lett.}, {\bf  12} (2005), no. 4,
  459--473.

\bibitem{Andersen-Kaneda2000} H. Andersen and M.
 Kaneda, On the $D$-affinity of the flag variety in type $B\sb 2$, {\em  Manuscripta Math.}, {\bf  103} (2000),
 no. 3, 393--399.

\bibitem{Am-Mac} M. F. Atiyah and I. G. Macdonald, {\em
Introduction to Commutative Algebra}, Addison-Wesley, Reading,
Mass., 1969.


\bibitem{JPAA-08invcharp} V. V. Bavula,   The inversion formulae for automorphisms of polynomial
algebras and differential operators in prime characteristic, {\em
J. Pure Appl. Algebra}, {\bf 212} (2008), 2320-2337.
(Arxiv:math.RA/0604477).

\bibitem{gendif} V. V. Bavula, Generators and defining relations for
the ring of differential operators on a smooth affine algebraic
variety,  (Arxiv:math.RA/0504475).

\bibitem{holmodp} V. V. Bavula, Dimension, multiplicity, holonomic modules, and an analogue of
the inequality of Bernstein for  rings of differential operators
in prime characteristic, (Arxiv:math. RA/0605073).


\bibitem{BGGDiffcone72} I. N. Bernstein, I. M. Gelfand and  S. I. Gelfand, Differential
 operators on a cubic cone. {\it Uspehi Mat. Nauk}  {\bf 27} (1972), no. 1, 185--190.

\bibitem{Bezrukavnikov-Mirkovic-Rumynin-AnnMAth} R. Bezrukavnikov, I. Mirkovic, and  D. Rumynin, Localization
of modules for a semisimple Lie algebra in prime characteristic,
 {\it Ann. of Math.}, to appear (Arxiv.org/math/0205144).


\bibitem{BogvadJA95} R. Bogvad, Some results on $\CD $-modules on
Borel varieties in characteristic $p>0$. {\it J. of Algebra} {\bf
173} (1995), 638--667.

\bibitem{BourbakiAlgII} N. Bourbaki, Algebra II. Chapters 4--7. Springer-Verlag, Berlin, 2003.


\bibitem{Eisenbook} D. Eisenbud, {\em Commutative algebra. With a
view toward algebraic geometry}. Graduate Texts in Mathematics,
150. Springer-Verlag, New York, 1995.

\bibitem{Fitting1936} H. Fitting, Die Determinantenideale eines
Moduls. {\em Jahrbericht der Deutschen Math.-Vereinigung}, {\bf
46} (1936), 195-229.

\bibitem{Haastert87} B. Haastert, \"{U}ber Differentialoperatoren und $D$-Moduln in positiver
Charakteristik. {\it Manuscripta Math.} {\bf 58} (1987), no. 4,
385--415.

\bibitem{Hun-Sharp93} C. L. Huneke, R. Y.  Sharp, Bass numbers of local cohomology modules.
 {\it Trans. Amer. Math. Soc.} {\bf  339} (1993), no. 2, 765--779.

\bibitem{Kashiwara-Lauritzen2002} M. Kashiwara and N.  Lauritzen,
Local cohomology and $\CD$-affinity in positive characteristic.
{\it  C. R. Math. Acad. Sci. Paris}, {\bf 335}  (2002), no. 12,
993--996.

\bibitem{LyubCrelle97} G. Lyubeznik,
 $F$-modules: applications to local cohomology and $D$-modules in characteristic $p>0$.
 {\it J. Reine Angew. Math.} {\bf  491} (1997), 65--130.


\bibitem{Ma} H. Matsumura, {\em Commutative ring theory.}
  Cambridge University Press,
Cambridge, 1994.


\bibitem{MR} J. C. McConnell and J. C. Robson, {\em Noncommutative Noetherian
rings}. With the cooperation of L. W. Small. Revised edition.
Graduate Studies in Mathematics, 30. American Mathematical
Society, Providence, RI, 2001.

\bibitem{Meb-Nar-MacLNM90} Z. Mebkhout and L. Narvaez-Macarro, Sur
les coefficients de Rham-Grothendieck des vari\'et\'es
alg\'ebriques, {\it Lecture Notes in Mathematics}, vol. 1454,
Springer-Verlag, Berlin/New York, 1990.


\bibitem{K.Smith95} K. E. Smith,  The $D$-module structure of $F$-split rings. {\it Math. Res. Lett.} {\bf  2}
 (1995), no. 4, 377--386.


\bibitem{KSm-vdB} K. E. Smith and M. van den Bergh, Simplicity
of rings of differential operators in prime characteristic. {\it
Proc. London Math. Soc.} {\bf 75} (1997), no. 1, 32--62.

\bibitem{SmStafDifopcurve} S. P. Smith  and J. T. Stafford, Differential
operators on an affine curve. {\it Proc. London Math. Soc.} {\bf
56} (1988), no. 2, 229--259.

\bibitem{Smith85LNM} S. P.  Smith, Differential operators on
commutative algebras, in {\em Ring Theory, LNM 1197}, (1985),
164-177.

\bibitem{Smith85LNM} S. P.  Smith, Differential operators on
 the affine and projective lines, in {\em Ring Theory, LNM 1220}, (1985), 157--177.



\end{thebibliography}
\end{document}